\documentclass[preprint,3p,fleqn,sort]{elsarticle}
\journal{Journal of Computational Physics}
\usepackage[labelfont=bf, singlelinecheck=false,
            justification=justified]{caption}
\usepackage{xcolor}
\usepackage{amsmath,amssymb}
\usepackage{subcaption}
\usepackage{mathrsfs}
\usepackage{stmaryrd} 
\usepackage{booktabs}
\usepackage{microtype}
\usepackage{mathtools} 
\usepackage{multirow}
\usepackage{enumitem}
\usepackage{setspace}
\usepackage{placeins} 
\usepackage{etoolbox,siunitx}
\usepackage{bm}
\robustify\bfseries
\makeatletter
\g@addto@macro\normalsize{%
	\setlength\abovedisplayskip{4pt}
	\setlength\belowdisplayskip{4pt}
	\setlength\abovedisplayshortskip{4pt}
	\setlength\belowdisplayshortskip{4pt}
}
\makeatother
\usepackage[T1]{fontenc}
\usepackage{accsupp} 
\SetSymbolFont{stmry}{bold}{U}{stmry}{m}{n}

\usepackage{lineno,hyperref}
\renewcommand{\vec}[1]{\ensuremath{\boldsymbol #1}}

\newcommand{\pderivative}[2]{\frac{\partial #1}{\partial #2}}
\newcommand{\mat}[1]{\ensuremath{\mathbf{#1}}}

\newcommand{\eg}{\textit{e.g.~}}

\newcommand{\etal}{\textit{et al.~}}
\renewcommand{\vec}[1]{\ensuremath{\boldsymbol #1}}
\newcommand{\He}{\mat{H}}
\newcommand{\vvec}{\vec{v}}

\newcommand{\jump}[1]{\ensuremath{\left\llbracket #1 \right\rrbracket}}

\newcommand{\halb}{\frac{1}{2}}

\let\rho\varrho

\begin{document}

\begin{frontmatter}

\title{Hybrid Entropy Stable HLL-Type Riemann Solvers for Hyperbolic Conservation Laws}

\author[aachen]{Birte Schmidtmann\corref{mycorrespondingauthor}}
\cortext[mycorrespondingauthor]{Corresponding author}
\ead{schmidtmann@mathcces.rwth-aachen.de}
\author[koeln]{Andrew R.~Winters}
\address[aachen]{MathCCES, RWTH Aachen University, Schinkelstr. 2, 52062 Aachen}
\address[koeln]{Mathematisches Institut, Universit\"at zu K\"oln, Weyertal 86-90, 50931 K\"oln}

\numberwithin{equation}{section}


\begin{abstract}	
	It is known that HLL-type schemes are more dissipative than schemes based on characteristic decompositions. However, HLL-type methods offer greater flexibility to large systems of hyperbolic conservation laws because the eigenstructure of the flux Jacobian is not needed. We demonstrate in the present work that several HLL-type Riemann solvers are provably entropy stable. Further, we provide convex combinations of standard dissipation terms to create hybrid HLL-type methods that have less dissipation while retaining entropy stability. The decrease in dissipation is demonstrated for the ideal MHD equations with a numerical example.
\end{abstract}

\end{frontmatter}
%
\section{Introduction}
%
We consider the numerical solution of systems of hyperbolic conservation laws of the form
\begin{equation}
	\label{eq:hypLaw}
	\pderivative{\vec{q}}{t} + \nabla\cdot\vec{f} = 0,
\end{equation}
on a domain $\Omega$. For a one-dimensional approximation we divide $\Omega$ into $K$ non-overlapping grid cells $C_i = [x_{i-1/2},\; x_{i+1/2}],\, i= 1,\ldots, K$ which are not necessarily equidistant. In the context of finite volume schemes, hyperbolic equations, such as \eqref{eq:hypLaw}, require a numerical flux function which fully determines the properties of the scheme \cite{LeVeque2002}. The numerical flux function takes as input the left and right value of $\vec{q}$ at the cell interface and solves a local Riemann problem. Smooth initial flows governed by \eqref{eq:hypLaw} may develop discontinuities (\eg shocks) in finite time. Thus, solutions are sought in the weak sense \cite{LeVeque2002}. 
However, weak solutions are not unique and need to be supplemented with additional admissibility criteria. Following the work of \eg \cite{Winters2016,Tadmor2003}, we use the concept of entropy to construct discretizations that agree with the second law of thermodynamics. That is, the numerical flux function will possess entropy stability, cf. \cite{Tadmor2003} and references therein.

In particular, we prove entropy stability for the HLL scheme and present the construction of HLL-type entropy stable numerical flux functions. It is known that HLL-type schemes are more dissipative than upwind schemes. However, HLL-type methods need less information about the eigendecomposition of the flux Jacobian. This is advantageous because the eigenstructure might be computationally expensive or no analytical expression exists, especially for large systems. As such, we consider three standard dissipation terms, namely Lax-Friedrichs (LF), HLL, and Lax-Wendroff (LW) and present two hybrid dissipation terms introduced in \cite{Schmidtmann2016}. We demonstrate that these five schemes are provably entropy stable. 

The paper is organized as follows: Sec.~\ref{scn:ESSolvers} provides a brief background on entropy stable numerical fluxes. In Sec.~\ref{scn:Riemann} we show entropy stability for the LF, HLL, and LW dissipation terms. The creation of two new hybrid entropy stable dissipation operators is shown in Sec.~\ref{scn:HybridRiemann}. We demonstrate in Sec.~\ref{sec:NumSim} that the new hybrid numerical flux reduce the overall dissipation in a standard finite volume scheme. Our conclusions and outlook are drawn in the final section.
%
\section{Entropy stable numerical flux functions}\label{scn:ESSolvers}
%
A numerical method that recovers the local changes in entropy as predicted by the continuous entropy conservation law is said to be \textbf{\emph{entropy conservative}}. Entropy conservation is only valid for smooth flow configurations. For discontinuous solutions, the entropy conservation law becomes the entropy inequality \cite{Tadmor2003}. A numerical scheme is said to be \textbf{\emph{entropy stable}} as long as the numerical approximation always obeys the entropy inequality
\begin{equation}
	\label{eq:EntropyInequality}
	\pderivative{S}{t} + \pderivative{{F}}{x} \le 0,
\end{equation}
where we assume that the system of hyperbolic conservation laws is equipped with a strongly convex mathematical entropy function, $S$, and a corresponding entropy flux, $F$, \cite{Tadmor2003}. 
%
It is known that without additional dissipation, entropy conservative numerical schemes produce high-frequency oscillations near shocks, see \eg \cite{Fjordholm2012,Winters2016}. Thus, for the approximation to remain valid for general flow configurations we must add a carefully designed dissipation term to ensure that \eqref{eq:EntropyInequality} discretely holds.

To create an entropy stable (ES) numerical approximation we start with a baseline entropy conserving (EC) numerical flux and then add a dissipation term. 
The resulting numerical flux at an arbitrary cell interface  $i+\halb$ takes the form
\begin{equation}
	\label{eq:entropystableflux}
	\vec{f}^{*,ES} = \vec{f}^{*,EC} - \frac{1}{2} \mat{D} \jump{\vec{q}},
\end{equation}
where $\vec{q}$ is the vector of conserved variables, $\mat{D} = \mat{D}(\vec q_{i}, \vec{q}_{i+1})$ is a suitable dissipation matrix evaluated at some mean state between the two cells, and  $\jump{\cdot} = (\cdot)_{i+1} - (\cdot)_{i}$ is the jump between the right and left cells. For simplicity of presentation we suppress the indices on the numerical flux, the dissipation matrix, and any jump terms. For an ES scheme, the baseline central flux from a classical Riemann solver is replaced by the baseline EC flux. To guarantee entropy stability, the dissipation term in \eqref{eq:entropystableflux} must be carefully designed to ensure that $\vec{f}^{*,ES}$ discretely satisfies the entropy inequality \eqref{eq:EntropyInequality}. To do so, we rewrite the dissipation term \cite{Roe2006}
\begin{equation}
	\label{eq:dissipation_term}
	\frac{1}{2} \mat{D} \jump{\vec{q}} \simeq \frac{1}{2} \mat{D} \He{} \jump{\vec{v}}, 
\end{equation}
with the vector of entropy variables $\vec{v}=\pderivative{S}{\vec{q}}$ and the entropy Jacobian $\He{} = \pderivative{\vec{q}}{\vec{v}}$ which relates the variables in conserved and entropy space. Substituting \eqref{eq:dissipation_term} into \eqref{eq:entropystableflux} the entropy stable numerical flux becomes
\begin{equation}
	\label{eq:entropystableflux2}
	\vec{f}^{*,ES} = \vec{f}^{*,EC} - \frac{1}{2} \mat{D}\He{} \jump{\vec{v}}.
\end{equation}
The reformulation of the dissipation term, incorporating the jump in entropy variables (rather than the jump in conservative variables) makes it possible to show entropy stability \cite{Barth1999}. From the structure of the entropy stable flux \eqref{eq:entropystableflux2}, we find a discrete version of the entropy inequality \eqref{eq:EntropyInequality} in cell $i$ to be \cite{Winters2016}
\begin{equation}
	\label{eq:EntropyInequalityDisc}
	\pderivative{S_i}{t} +  \left( F_{i+\halb} - F_{i-\halb} \right) \le -\halb\jump{\vvec}^T\mat{D}\He{}\jump{\vvec}\stackrel{!}{\leq} 0.
\end{equation}
Thus, to guarantee discrete entropy stability, it is sufficient to show that $\mat{D}\mat{H}$ is symmetric positive definite (s.p.d). 
%
%
\section{Entropy stable classical Riemann solvers}\label{scn:Riemann}
%
In this section we demonstrate entropy stability for the numerical flux of the form \eqref{eq:entropystableflux2} for the dissipation matrix $\mat{D}$ of the LF, HLL, and LW scheme. To do so, we first assume that the flux Jacobian, $\mat{A}$, or a suitable Roe matrix, exists with the properties
\begin{equation}
	\label{eq:assume}
	\mat{A} = \mat{R}\boldsymbol\Lambda\mat{R}^{-1},\quad 
	\He{} = \left(\mat{R}\mat{Z}\right)\left(\mat{R}\mat{Z}	\right)^{T},
\end{equation}
where $\mat{R}$ is the eigenvector matrix, $\boldsymbol\Lambda$ the diagonal corresponding eigenvalue matrix, and $\mat{Z}$ is a positive diagonal scaling matrix which creates a set of entropy scaled eigenvectors $\mat{R}\mat{Z}$ \cite{Barth1999}. We see that, by construction in \eqref{eq:assume}, the matrix $\He{}$ is s.p.d.
In the later proofs we only use the existence of the matrices $\mat{R}$, $\boldsymbol\Lambda$, and $\mat{Z}$, whereas, in practice, their explicit form does not need to be known. This is advantageous, because for large systems of conservation laws, the eigendecomposition is expensive to compute or is not available.

To write the entropy stable LF scheme we substitute the dissipation matrix 
\begin{align}
	\mat{D}_\text{LF}=\frac{\Delta x}{\Delta t}\,\mat{I},
\end{align}
into the form \eqref{eq:entropystableflux2}.
Now the complete dissipation term for \eqref{eq:entropystableflux2} only depends of the known s.p.d matrix $\mat{H}$, so discrete entropy stability \eqref{eq:EntropyInequalityDisc} for $\mat{D}_{\text{LF}}$ follows immediately.
We note that under the same assumption the local Lax-Friedrichs (LLF) and Roe-type dissipation terms satisfy the discrete entropy stability \cite{Winters2016}.
Next, we consider the HLL flux. The numerical flux function of ES-HLL can be written in form \eqref{eq:entropystableflux2}, with the dissipation matrix
\begin{equation}\label{eq:DHLL}
	\mat{D}_\text{HLL} = a_0\mat{I} + a_1\mat{A}, \quad a_0 = \frac{|\lambda_L|\lambda_R - |\lambda_R|\lambda_L}{\lambda_R-	\lambda_L}, \quad a_1 = \frac{|\lambda_R| - |\lambda_L|}{\lambda_R-\lambda_L}.
\end{equation}
Here, $\lambda_{L,R}$ are the fastest signal velocities with $\lambda_L < \lambda_R$. 
>From assumption \eqref{eq:assume} it is straightforward to show that the discrete entropy stability condition \eqref{eq:EntropyInequalityDisc} is equivalent to showing that $a_0 +a_1\lambda_i\geq 0$ for all $\lambda_i\in[\lambda_L,\lambda_R]$. From the form of the coefficients $a_0$ and $a_1$, keeping in mind that $\lambda_L < \lambda_R$ and $\lambda_i\in[\lambda_L,\lambda_R]$, we find
\begin{equation}
	\label{eq:HLLcond}
		|\lambda_L|\lambda_R - |\lambda_R|\lambda_L + \left( |\lambda_R| - |\lambda_L|\right)\lambda_i\geq 0\qquad
		\Leftrightarrow \qquad |\lambda_L|\left(\lambda_R - \lambda_i\right) + |\lambda_R|\left(\lambda_i - \lambda_L\right)\geq 0.
\end{equation}
Finally, we consider the dissipation matrix for LW
\begin{equation}
	\label{eq:DLW}
	\mat{D}_\text{LW}= \frac{\Delta t}{\Delta x}\mat{A}^2,
\end{equation}
in the ES numerical flux \eqref{eq:entropystableflux2} and find that the discrete entropy inequality \eqref{eq:EntropyInequalityDisc} holds, i.e.,
\begin{equation}
	\label{eq:LWcond}
	\pderivative{S}{t} +  \left( F_{i+\halb} - F_{i-\halb} \right)
	\leq -\halb \frac{\Delta t}{\Delta x}\left(\mat{Z}\mat{R}^T \jump{\vvec}\right)^T\boldsymbol\Lambda^2 \left(\mat{Z}\mat{R}^T \jump{\vvec}\right)
	\leq 0.
\end{equation}
%
\section{Hybrid entropy stable Riemann solvers}\label{scn:HybridRiemann}
%
In this section, we consider dissipation matrices for hybrid Riemann solvers constructed in \cite{Schmidtmann2016}, motivated by the work of Degond \etal \cite{DegondPeyrardRussoVilledieu1999}. The hybrid solvers were constructed using weighted combinations of the dissipation matrices described in Sec.~\ref{scn:Riemann}. The weighting is chosen to reduce dissipation, especially for signal velocities close to zero. As such, the hybrid terms contain a parameter $\omega \in [0,1]$ which allows further control over the amount of dissipation added to the scheme.
%
\subsection{$HLL\omega$}\label{scn:HLLohm}
%
First we consider a generalization of the ES-HLL flux. The dissipation matrix, containing a parameter $\omega\in[0,1]$, is a linear function of the eigenvalues of the associated flux Jacobian (or Roe matrix) $\mat{A}$ 
\begin{equation}\label{eq:DHLLw}
	\begin{aligned}
		&\mat{D}_{HLL\omega} = b_0(\omega)\mat{I} + b_1(\omega)\mat{A},\\
		& b_0(\omega) = \dfrac{\lambda_R( \omega\,\lambda_L^2+(1-\omega) |\lambda_L|)-\lambda_L( \omega\, \lambda_R^2+(1-\omega) |\lambda_R|)}{\lambda_R-\lambda_L}, \quad
				b_1(\omega) = \dfrac{(1-\omega)(|\lambda_R|-|\lambda_L|)+\omega(\lambda_R^2-\lambda_L^2)}{\lambda_R-\lambda_L}.
	\end{aligned}
\end{equation}
We note that if $\omega=0$ we recover the ES-HLL matrix \eqref{eq:DHLL}. For $\mat{D}_{\text{HLL}\omega}$ to fulfill the discrete entropy stability condition \eqref{eq:EntropyInequalityDisc}, it is sufficient to show that $b_0(\omega) + b_1(\omega) \lambda_i\geq 0$ holds for all $\lambda_i\in[\lambda_L, \lambda_R]$. 
Inserting the coefficients $b_0(\omega),\, b_1(\omega)$, this condition rearranges to become
\begin{equation}\label{eq:HLLwcond}
	\omega \left[ \lambda_L^2\left(\lambda_R - \lambda_i\right) + \lambda_R^2\left(\lambda_i - \lambda_L\right)\vphantom{\frac{.}{.}} \right] + 
	(1-\omega) \left[\vphantom{\frac{.}{.}}|\lambda_L|\left(\lambda_R - \lambda_i\right) + |\lambda_R|\left(\lambda_i - \lambda_L\right)\right] \geq 0,\quad\forall\, \omega \in[0,1].
\end{equation}
%
\subsection{$HLLX\omega$}\label{scn:HLLXohm}
%
\begin{figure}
	\centering
	\includegraphics[width=0.5\textwidth]{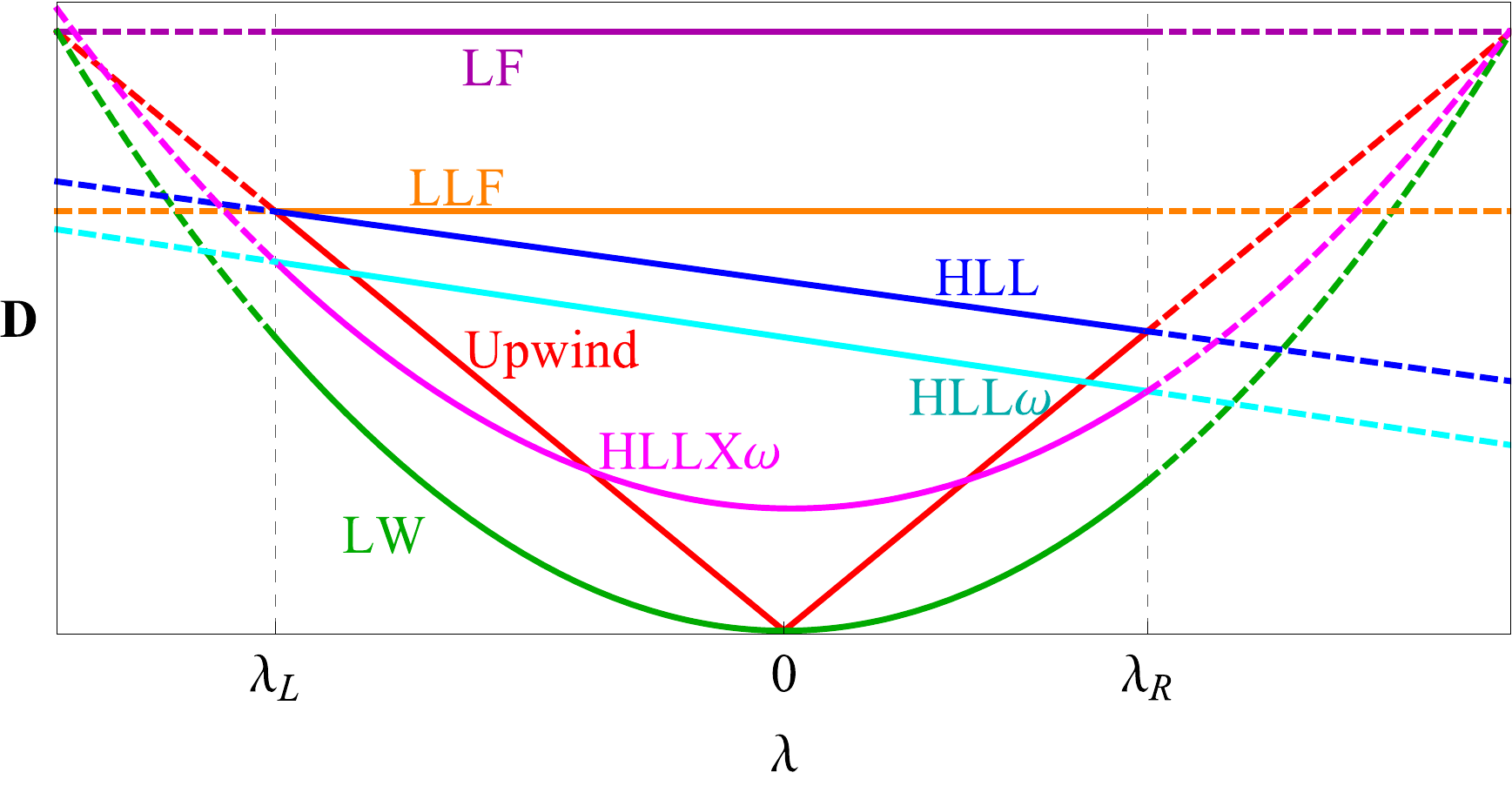}
	\caption{Visualization of the dissipation matrix $\mat{D}$ as a function of the eigenvalue $\lambda$. For the hybrid dissipation matrices we select $\omega=0.4$.}
	\label{fig:Riemann}
\end{figure}
Next, we treat a hybrid dissipation matrix, HLLX$\omega$, that includes the parameter $\omega$ and the quadratic LW term. This term requires squaring the flux Jacobian (or applying the flux twice) but reduces the overall magnitude of the dissipation, see Fig.~\ref{fig:Riemann}. 
The dissipation matrix of HLLX$\omega$ is a weighted combination of LF, HLL$\omega$, and LW, see \cite{Schmidtmann2016} for more details.
\begin{equation}
	\begin{aligned}	
		&\mat{D}_{HLLX\omega} = \beta_0(\omega) \mat{D}_{LF} + \beta_1(\omega) \mat{D}_{HLL\omega} + \beta_2(\omega) \mat{D}_{LW},\\
		&\beta_0(\omega) = \beta(\omega)\frac{(1-\omega)\left|\lambda_L\,\lambda_R\right|}{(1-\omega) + \omega  (\left|\lambda_L\right|+\left|\lambda_R\right|)},\quad
		\beta_1(\omega) =1 - \beta(\omega)\left(\frac{1-\omega}{\left|\lambda_L\right|+\left|\lambda_R\right|}+\omega \right)^{-1},\\
		&\beta_2(\omega) =\beta(\omega),\quad
		\beta(\omega) = \omega + (1-\omega)\frac{\lambda_R-\lambda_L-\big| |\lambda_R| - |\lambda_L| \big|}{(\lambda_R-\lambda_L)^2}\equiv \omega + (1-\omega)\alpha.
		\label{eq:beta}	
	\end{aligned}	
\end{equation}
The resulting numerical flux is not strictly monotone. However, we can guarantee entropy stability. We have already shown discrete entropy stability for LF, HLL$\omega$ and LW. Thus, to demonstrate ES for HLLX$\omega$, it suffices to show the positivity of the coefficients $\beta_i(\omega), \,i=0,1,2$. Note that for $\omega=1$ we obtain LW, which has already been shown to be discretely ES. Thus, let us consider $\omega\in [0,1)$. The claim holds since
\begin{equation}
	\beta_2(\omega) = \beta(\omega) \geq 0 \quad \Leftrightarrow\quad \omega + (1-\omega)\alpha\geq 0\quad  \Leftrightarrow\quad \alpha \geq 0\quad\Leftrightarrow\quad \lambda_R-\lambda_L-\big| |\lambda_R|-|\lambda_L|\big|\geq 0.
\end{equation}
Hence, we directly see that $\beta_0\geq 0$ because it is a combination of non-negative terms. In order to show that $\beta_1(\omega)\geq 0$ we consider the equivalent expression
\begin{equation}
	\beta(\omega)\leq \left(\frac{1-\omega}{\left|\lambda_L\right|+\left|\lambda_R\right|}+\omega \right)\quad\Leftrightarrow\quad \omega + (1-\omega)\alpha \leq \left(\frac{1-\omega}{\left|\lambda_L\right|+\left|\lambda_R\right|}+\omega \right)\quad\Leftrightarrow\quad \alpha (\left|\lambda_L\right|+\left|\lambda_R\right|)\leq 1.
\end{equation}
We distinguish two cases. Case 1: $\lambda_L$ and $\lambda_R$ are both positive or both negative. Then $\lambda_R-\lambda_L=\big| |\lambda_R|-|\lambda_L|\big|$ and therefore $\alpha=0$. Case 2: If $\lambda_L$ and $\lambda_R$ are of opposite sign, then $\lambda_R-\lambda_L = |\lambda_R|+|\lambda_L|$, which yields
\begin{equation}
\alpha\left( |\lambda_R|+|\lambda_L|\right) = 1- \big||\lambda_R|-|\lambda_L|\big|\left( |\lambda_R|+|\lambda_L|\right)/\left(\lambda_R-\lambda_L\right)^2 
	\leq 1.
\end{equation}
%
\section{Numerical example - application to ideal MHD}\label{sec:NumSim}
%
We now apply the new hybrid ES-HLLX$\omega$ Riemann solver to the equations of ideal magnetohydrodynamics (MHD). The ideal MHD equations are a hyperbolic system that decribes the flow of plasma assuming infinite electric resistivity, see \eg \cite{Winters2016}.
As a proof of concept we implement the hybrid ES-HLLX$\omega$ solver into a first order finite volume framework. We use a 1D shock tube problem for the ideal MHD equations to demonstrate the reduced dissipation of the new hybrid numerical flux. We consider the magnetic shock tube of Torrilhon \cite{Torrilhon2003}
\begin{equation}\label{eq:magShock}
	\begin{bmatrix}
		\rho, \rho u, \rho v, \rho w, p, B_1, B_2, B_3 
	\end{bmatrix}^T 
	= \left\{
	\begin{array}{lc}
		\left[1,0,0,0,1,1.5,0.5,0.6\right]^T, & \textrm{if}\quad x \leq 0, \\ 
		\left[1,0,0,0,1,1.5,1.6,0.2\right]^T, & \textrm{if}\quad x > 0,
	\end{array}
	\right.
\end{equation}
on the domain $\Omega = [-4,4]$ with Dirichlet boundary conditions and an adiabatic index $\gamma=\frac{5}{3}$. The baseline EC flux needed for the numerical flux ansatz \eqref{eq:entropystableflux} is chosen to be the entropy conserving and kinetic energy preserving flux found in \cite[App. B]{Winters2016}. Briefly, we note that the non-linear stability of the scheme depends on which underlying baseline EC flux is chosen to build the scheme. However, in our experience, the non-linear stability properties of a scheme created with the available EC baseline fluxes is nearly identical for low Mach number test cases (like the magnetic shock tube). We compare the ES-Roe flux of \cite{Winters2016} and ES-HLLX$\omega$ with $\omega=0.925$. 

\begin{figure}[t]
	\begin{center}
		{
		\includegraphics[scale=0.55]{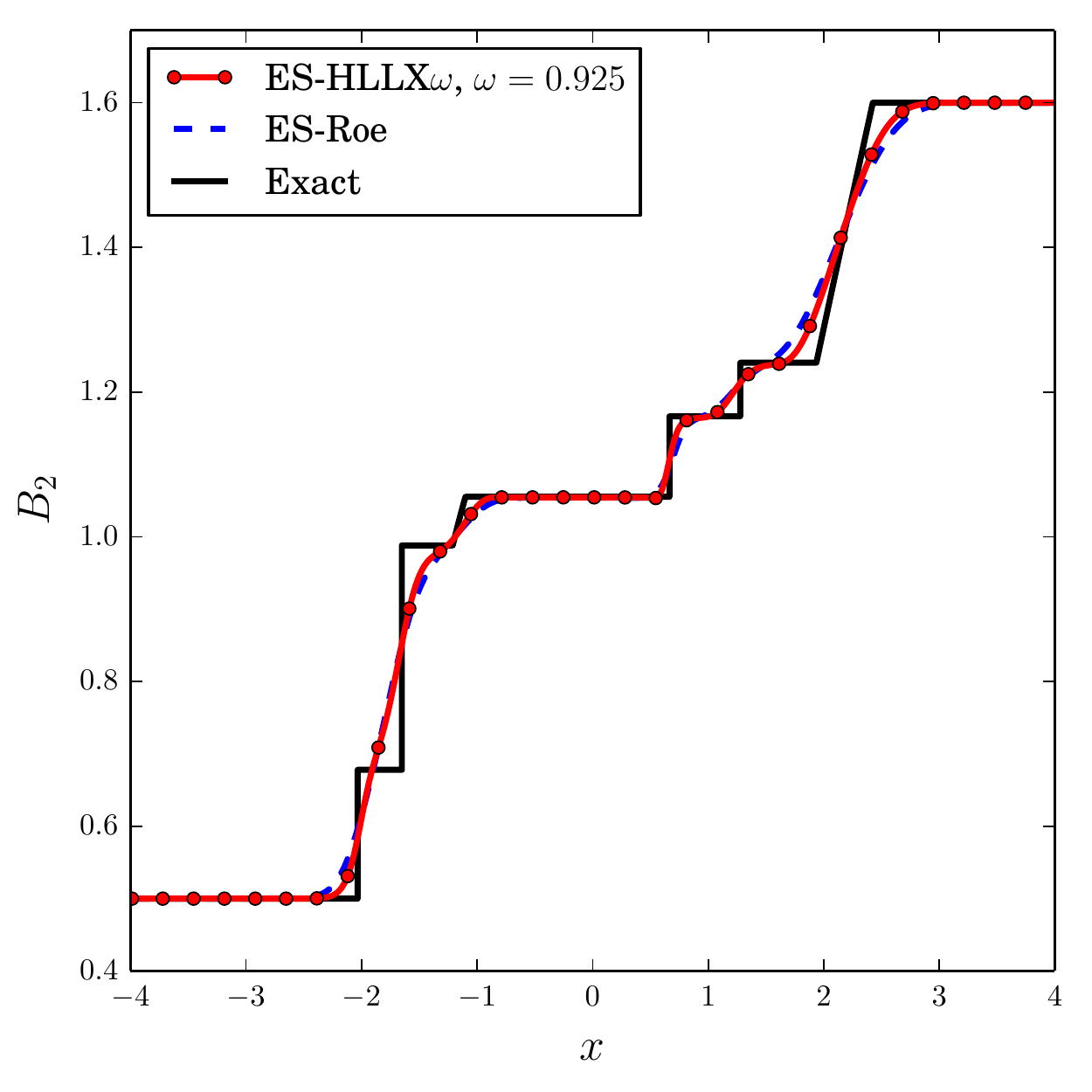}
		}
		{
		\includegraphics[scale=0.55]{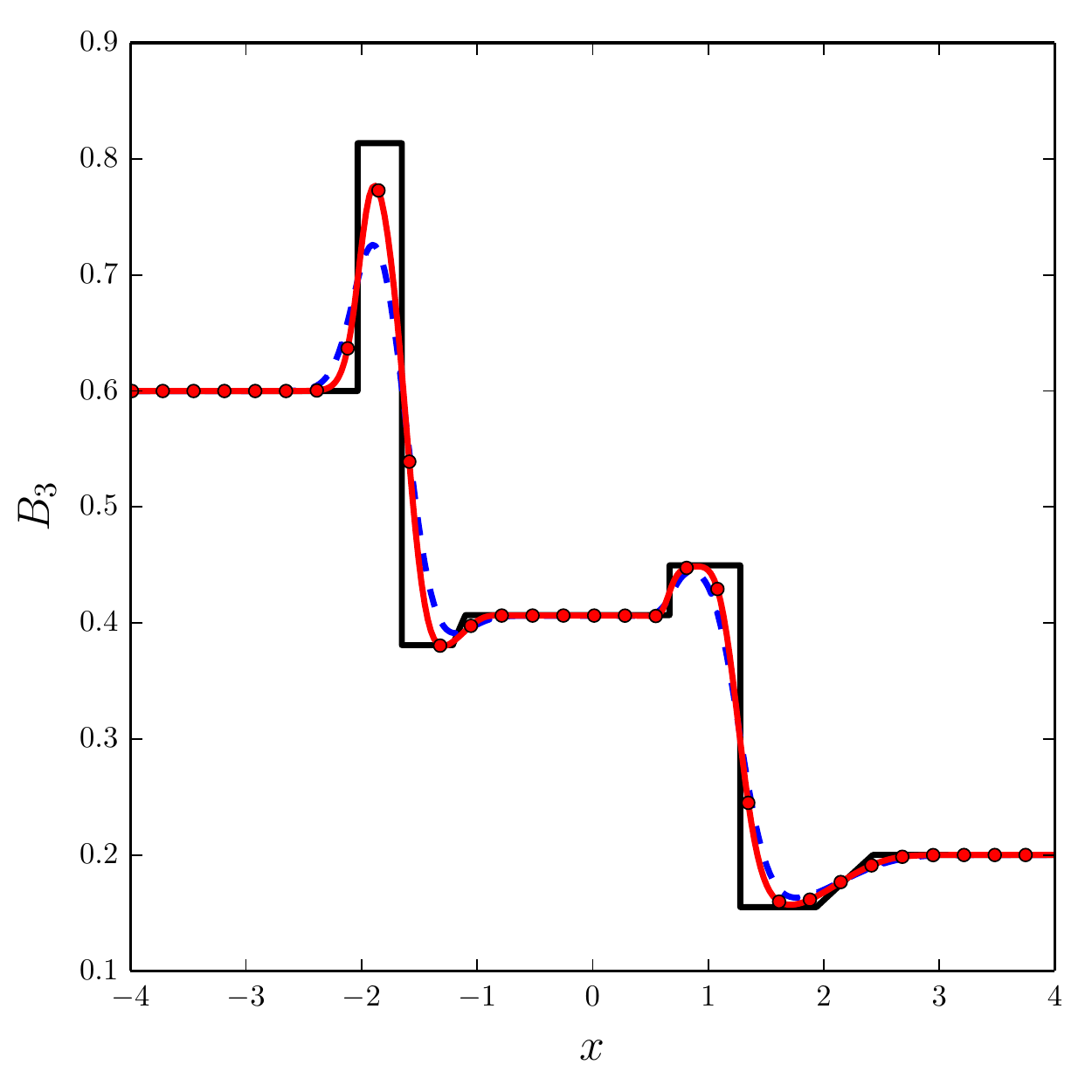}
		}
		\caption{Comparison of the computed solution of $B_2$ and $B_3$ using ES-Roe (dashed) and ES-HLLX$\omega$ (solid with knots) with $\omega = 0.25$ for the magnetic shock tube problem \eqref{eq:magShock} at $T=1.0$ on 300 regular grid cells.}
		\label{fig:MagTube}
	\end{center}
\end{figure}

Each dissipation term $\mat D = \mat D(\vec{q}_{i}, \vec{q}_{i+1})$ in the ES-HLLX$\omega$ dissipation matrix \eqref{eq:beta} is created using a simple arithmetic mean state of the primitive variables. We note that the value of $\omega$ could be chosen adaptively using for example a pressure switch \cite{Chandrashekar2012}. 
Fig.~\ref{fig:MagTube} presents the computed solution of $B_2$ and $B_3$ on 300 regular grid cells against the exact solution of the Riemann problem at the final time $T=1.0$. We see that the entropy stable hybrid numerical flux has less dissipation than the entropy stable Roe-type scheme. 
%
\section{Conclusion}\label{scn:Conclusion}
%
In this work we constructed two one-parameter families of hybrid entropy stable numerical fluxes. An advantage of the new numerical flux functions is that they remain applicable even when the eigenstructure of the flux Jacobian matrix is unknown. 
The derivations and proofs in this work are kept general such that the hybrid entropy stable solvers can be applied to a broad range of non-linear hyperbolic conservation laws. As an example, we applied the novel numerical fluxes to the ideal MHD equations and demonstrated the decreased magnitude of dissipation for the hybrid solvers versus a standard solver. 
In the future we plan to apply the hybrid entropy stable Riemann solvers to other complex hyperbolic systems, such as the two-layer shallow water or the regularized 13-Moment Equations of Grad.

%

\section*{References}
\bibliography{./mybibfile}

\end{document}